\documentclass[12pt,a4paper]{article}
\usepackage[cp1251]{inputenc}
\usepackage{graphicx}
\usepackage{amssymb}
\usepackage{amsmath}
\usepackage{amsthm}
\usepackage{euscript}
\numberwithin{equation}{section}
\newtheorem{theorem}{Theorem}
\newtheorem{lemma}{Lemma}

\newtheorem{rem}{Remark}

\newtheorem{ex}{Example}

\newcommand{\R}[0]{\mathbb{R}}

\newcommand{\Z}[0]{\mathbb{Z}}

\newcommand{\be}[0]{\begin{equation}}
\newcommand{\ee}[0]{\end{equation}}
\newcommand{\bez}[0]{\begin{equation*}}
\newcommand{\eez}[0]{\end{equation*}}
\newcommand{\bl}[0]{\begin{lemma}}
\newcommand{\el}[0]{\end{lemma}}

\newcommand{\ep}[0]{$\hspace{\fill} \square$}
\newcommand{\paragraf}[1]{\par
\bigskip{\centerline{\bf #1}}\medskip}

\newcommand{\abs}[1]{\begin{quotation} {\small
 \centerline{{\bf Abstract}}  \smallskip 
#1}
\end{quotation}}
\author{I.\,N.~Shnurnikov\footnote{NRU HSE, Moscow}}

\title{On the number of connected components of complements to arrangements of subtori}
\date{}
\begin{document}
\maketitle
\abs{
\noindent
 We consider the arrangements of subtori in a flat $d$--dimensional torus $T^d$. Let us consider an arrangement on $n$ subtori of codimension one, let $f$ be the number of connected components of the complement in $T^d$ to the union of subtori. We found the set of all possible numbers $f$ for given $n$ and $d$ and arbitrary arrangements of subtori.
 }
 
{\bf MSC:} 52C35

{\bf Keywords:} arrangements of subtori, number of connected components.

\paragraf{Introduction}
The theory of plane arrangements in affine or projective spaces has been investigated thoroughly, see the book of Orlic and Terao, \cite{Orlic_1992} and Vassiliev's review \cite{Vasil'ev}. Inspired by the conjecture of Gr\"unbaum \cite{Grunbaum 72}, Martinov \cite{Martinov 93} found all possible pairs $(n,f)$ such that there is a real projective plane arrangement of $n$ pseudolines and $f$ regions. It turns out, that some facts concerning the arrangements of hyperplanes or oriented matroids could be generalized to arrangements of submanifolds, see Deshpande \cite{Deshpande}.

We study the sets $F(T^d,n)$ of connected components numbers of the complements in a flat $d$--dimensional torus $T^d$ to the unions of $n$ closed connected codimensional one subtori. Author in \cite{Shnurnikov 11} found the sets $F(T^2,n)$ of region numbers in arrangements of $n$ closed geodesics in a two dimensional torus and a Klein bottle with locally flat metrics.
The main result of the present paper is the following:
$$
F\left(T^d,n\right) = \{l \in \mathbb{N} \quad | \quad n-d+1 \leq l \leq n \ \ \text{or} \  \ 2(n-d)\leq l\}
$$
for $n>d \geq 2$. Also $F\left(T^d,n\right)=\mathbb{N}$ for $2\leq n \leq d$.  $F\left(T^d,1\right)=\{1\}$.

\paragraf{Main part.}

Let $\R^d$ be the real $d$--dimensional vector space and $u_1, \dots, u_d$ are linear independant vectors in $\R^d$. Let $L$ be the set of all linear combinations $k_1u_1+ \dots+k_du_d$ with integer coefficients $k_1, \dots, k_d$. The set $L$ could be considered as the group of  transitions, acting on $\R^d$.
The quotient space of this action is torus, the induced metric (from euclidean metric in $R^n$) is locally flat. So we call this quotient space {\it flat torus $T^d$} and denote by $\R^d/L$. So the factorisation map $\varphi: \R^d \to \R^d/L $ sends point $x$ to equivalent class $x+L$.  
Let $(x_1, \dots, x_d)$ are coordinates of $\R^d$ accordying to the basis $u_i$. Let us consider a hyperplane $K$ in $\R^d$ with equation $\sum a_ix_i =c$. The subset $\varphi(K) \subset T^d$ is closed iff the numbers $\frac{a_i}{a_j}$ are rational for all pairs $i,j$ such that $a_j\neq 0$. In the set $L$ one may choose another basis $v_1, \dots, v_d$, then the matrix of basis transformation is a $SL(\Z)$ matrix. 

Let $A$ be the union of $n$ codimensional one subtori in the flat $d$--dimensional torus $T^d$. Consider the connected components of the complement $T^d\setminus A$; denote the number of connected components by $f=|\pi_0\left(T^d\setminus A\right)|$; let $F\left(T^d,n\right)$ be the set of all possible numbers $f$.


\begin{lemma} Let  $K$ be a hyperplane in $\R^d$ with equation $\sum a_ix_i =c$ in a basis $U=(u_i)$ of the lattice  for integers $a_i$ such that $\gcd(a_1, \dots, a_d)=1$. 

(a) There is a new basis $V=(v_i)$ of the lattice, accordying to which $K$ is given by equation $y_1=0.$ The matrix $M$, $V=UM$ of basis transfer is a $SL(\Z)$ matrix with integer coefficients and the determinant is 1.

(b) Let $L$ be such that $U$ is an orthogonal basis and $|u_i|=1$ for all $i$. Then the distance between $K$ and the nearest point $x \in L$ of the lattice {\small }such that $x \notin K$ is $\frac 1{\sqrt{\sum_{i=1}^d a^2_i}}$.

(c) Let $U$ be an orthogonal basis and $|u_i|=1$ for all $i$. The volume of $T^d$ is $1$. The subset $\varphi (K)$ is a subtorus in $T^d$ and the $(d-1)$--dimensional volume of $\varphi(K)$ is $\sqrt{\sum_{i=1}^d a^2_i}$.
\end{lemma}

\begin{lemma} 
\label{f_via_nm}
Let $A_1, \dots A_n$ be the set of $n$ codimensional one subtori in the flat $d$--dimensional torus $T^d$. Let $m$ be the maximal number of parallel tori among $A_1, \dots A_n$. Then $f\geq m(n-m-d+2)$. 
\end{lemma} 
\proof Induction on $d$, base $d=1$ is trivial. Suppose the statement is true for $d-1$ and every $n$. Let us prove the statement for $d$. Consider $m$ parallel subtori, they divide torus into $m$ components $U_i$, $i=1, \dots, m$ and each component is homeomorphic to prime product of a segment and $(d-1)$--dimensional torus. For each $U_i$ let us consider a control $(d-1)$--dimensional torus $R_i$ such that the other $n-d$ tori of $A_1, \dots A_n$ form on $R_i$ an arrangement of $n-m$ subtori and so by the induction assumption they divide $R_i$ into at least $n-m-d+2$ components. So, each $U_i$ is divided by the other $n-d$ tori of $A_1, \dots A_n$ into at least $n-m-d+2$ regions and the induction is over.
\ep

\begin{lemma} 
\label{main_for_torus}
Let $n>d \geq 2.$ Then $f \geq 2n-2d $ or: $f \leq n$ and there are at least  $n-d+1$ parallel tori.
\end{lemma}
\proof Induction on $d$, base $d=2$ is found in \cite{Shnurnikov 11}. Suppose the statement is true for $d-1$, let us prove it for $d$. Let $m$ be the maximal number of parallel tori in the arrangement $A=\{A_1, \dots, A_n\}$. By lemma \ref{f_via_nm} 
we have $f \geq m(n-m-d+2)$. So we over in case $2 \leq m \leq n-d$, because in this case  $m(n-m-d+2) \geq 2(n-d)$. 

Let us consider the case $m=1$. Then we consider a control torus $R_1$ which is parallel to torus $A_1$ and is in general position to other tori of arrangement. Let us consider arangement $B=\{B_1, \dots, B_{n-1}\}$, which is formed in $R_1$ by other tori of arrangement $A$. We know that $f(A)\geq f(B)$. If in $B$ there are at most $n-d$ parallel tori then by induction assumption we have $f(B) \geq 2(n-d)$. Otherwise in $B$ there is at least $n-d+1$ parallel tori and this means that there exist a $(d-2)$--dimensinal torus $S^{d-2} $which is parallel to at least $n-d+2$ tori of $A$. Consider the two--dimensional subtori $W^2$ which is orthogonal to $S^{d-2}$. Then by induction assumption for $d=2$ we have that $A$ divides $W^2$ into $f(W^2)\geq 2(n-d)$ regions and so  $f(A)\geq f(W^2) \geq 2(n-d)$.

Let us consider the case $m\geq n-d+1$. The subset of $m$ parallel tori divide $T^d$ into $m$ regions. If at least one of them is divided into more then one connected components by the other tori of $A$, then every of $m$ regions is divided in at least two connected components and so $f(A)\geq 2m \geq 2(n-d)$. In the other case the induction is over.
\ep 	

\begin{theorem}
\label{theorem F set for d-dimens. torus}
$$
F\left(T^d,n\right) = \{l \in \mathbb{N} \quad | \quad n-d+1 \leq l \leq n \ \ \text{or} \  \ 2(n-d)\leq l\}
$$
for $n>d \geq 2$. Also $F\left(T^d,n\right)=\mathbb{N}$ for $2\leq n \leq d$.  $F\left(T^d,1\right)=\{1\}$.
\end{theorem}
\proof 
The case $d=2$ is proved in \cite{Shnurnikov 11}. Now we assume $d \geq 3$. 
From the lemma \ref{main_for_torus}, lemma \ref{f_via_nm} we see that the set $F\left(T^d,n\right)$ does not contain other integers.

We construct examples for $\leq n$ and $\geq 2n-2d$ regions separately.

Let us consider $n$ hyperplanes  in $\R^d$ (an equation corresponds to a hyperplane):
\begin{gather*}
x_i=0, \quad 1 \leq i \leq k,\\
x_{k+1}=c_{i-k}, \quad k+1 \leq i \leq n
\end{gather*}
for some integer $k$, $0 \leq k \leq d-1$ and real $c_{i-k}$ with different fractional parts. By the factorization map $\R^d \to \R^d/Z^d$ we get a set  $\left\{T_i^{d-1}, i=1,\dots,n\right\}$ of $n$ codimensional one subtori. And the complement is homeomorphic to the prime product
$$
T^d\setminus\bigcup_i T_i^{d-1}\approx \R^k\times\left(S^1\setminus\{p_1,\dots,p_{n-k}\}\right)\times\left(S^1\right)^{d-k-1},
$$
where $S^1\setminus\{p_1,\dots,p_{n-k}\}$ denotes a circle without $n-k$ points. Hence the number of complement regions equals $n-k$, for an integer $k$ such that $0 \leq k \leq d-1$.

Now let us take a nonnegative integer $k$ and construct an arrangement with $2n-2d+k$ connected components of the complement. We determine the subtori by equations:
\begin{gather*}
x_i=0, \quad \text{for} \quad 2\leq i \leq d,\\
x_2=kx_1+ \frac12,\\
x_1=c_j \quad \text{for} \quad j=1, \dots, n-d,
\end{gather*}
whereas numbers $kc_j+\frac 12$ are not integer for any $j$. ( This means that the intersection of three subtori
$$
x_2=kx_1+ \frac12, \ x_1=c_j, \ x_2=0
$$
is an empty set.) Therefore
$$
T^d\setminus\bigcup_{i=3}^d\{x_i=0\}\approx T^2\times \R^{d-2}.
$$
In the 2--dimensional torus the equations
\begin{gather*}
x_2=0,\\
x_2=kx_1+\frac 12,\\
x_1=c_j \ \text{for} \ j=1, \dots, n-d
\end{gather*}
produce the arrangement of $n-d+2$ closed geodesics. The geodesic union divides the torus into $2n-2d+k$ connected components (for more details on the arrangements of closed geodesics in a flat torus see \cite{Shnurnikov 11}).
\ep

\begin{lemma}
Let $K_1$ and $K_2$  be hyperplanes in $\R^d$ given by the equations $x_1=0$ and $\sum_i a_ix_i=c$ where $a_i$ are integers. Then the intersection of subtori $\varphi(K_1)$ and $\varphi(K_2)$ in $T^d$ consists of $\gcd(a_2, \dots, a_d)$ connected components, each of which is a $(d-2)$--dimensional subtorus.
\end{lemma}

\begin{lemma}
Let $\mathbf{a}=(a_1, \dots, a_n)$ and $\mathbf{b}=(b_1, \dots,b_n)$ be integer vectors such that $\gcd(a_1, \dots, a_n)=1$ and $\gcd(b_1, \dots, b_n)=1$. Let two subtori in a flat $d$--dimensional torus are given by vectors $\mathbf{a}$ and $b$. Let $f(a,b)$ be the number of connected components of the intersection of two subtori. Then 
\begin{multline*}
f(\mathbf{a},\mathbf{b})=\gcd \Bigg(  \frac{a_2b_1-b_2a_1}{\gcd(a_1,a_2)}, \frac{a_3\big(b_1u_1^{(1)}+b_2u_2^{(1)}\big)-b_3\gcd(a_1,a_2)}{\gcd(a_1,a_2,a_3)}, \\
\dots \\
\frac{a_n\big(b_1u_1^{(n-2)}+b_2u_2^{(n-2)}+\dots +b_{n-1}u_{n-1}^{(n-2)}\big)-b_n\gcd(a_1, \dots, a_{n-1})}{\gcd(a_1, \dots, a_n)} \Bigg)
\end{multline*}
where and $u_i^j$ are integers such that 
$$
\gcd(a_1, \dots, a_{j+1})=a_1u_1^{(j)}+a_2u_2^{(j)}+ \dots + a_{j+1}u_{j+1}^{(j)}
$$
for each $j=1, \dots, n-2$.
\end{lemma}


\begin{thebibliography}{99}

\bibitem{Deshpande}
P.~Deshpande, Arrangements of Submanifolds and the Tangent Bundle Complement. {\it Electronic Thesis and Dissertation Repository}, Paper 154 (2011).

\bibitem{Grunbaum 72}
B.~Gr$\mathrm{\ddot{u}}$nbaum, {\it  Arrangements and Spreads.} AMS, Providence, Rhode Island, 1972.

\bibitem{Martinov 93}
N.~Martinov, Classification of arrangements by the number of their cells. {\it Discrete and Comput. Geometry} (1993) {\bf 9}, N 1, 39--46.


\bibitem{Orlic_1992}
P.~Orlic, H.~Terao, Arrangements of Hyperplanes. Springer -- Verlag, Berlin -- Heidelberg, 1992. 329 pp.

\bibitem{Shannon 76}
R.W.~Shannon, A lower bound on the number of cells in arrangements of hyperplanes. {\it Jour. of combinatorial theory (A)}, {\bf 20}, (1976) 327--335.


\bibitem{Shnurnikov 10}
I.\,N.~Shnurnikov, Into how many regions do $n$ lines divide the plane if at most $n - k$ of them are concurrent?  {\it Moscow Univ. Math. Bull., ser. 1} (2010) {\bf 65}:5, 208~--~212.


\bibitem{Shnurnikov 11}
I.\,N.~Shnurnikov, On the number of regions formed by arrangements of closed geodesics on flat surfaces,
{\it Math. Notes} {\bf 90}, N 3 -- 4 (2011), 619 -- 622.

\bibitem{Vasil'ev}
V.\,A.~Vassiliev, Topology of plane arrangements and their complements. {\it Uspekhi Mat. Nauk},  {\bf 56}, iss. 2(338), (2001), 167 --- 203.


\end{thebibliography}
\end{document}